\documentclass{article}
\usepackage{amsmath,amsfonts,verbatim,amssymb,fullpage, amsthm}
\usepackage{color,graphics}
\usepackage{enumerate}
 \DeclareMathOperator{\diam}{diam}

 \DeclareMathOperator{\supp}{supp}

\newcommand\remove[1]{}

\newcommand{\E}{\mathbb{E}}

\newcommand{\R}{\mathbb{R}}

\newcommand{\e}{\varepsilon}

\newcommand{\jnote}[1]{\marginpar{J}$\ll$\textbf{#1 --James}$\gg$}

\theoremstyle{plain}
\newtheorem{lemma}{Lemma}[section]
\newtheorem{theorem}[lemma]{Theorem}

\newtheorem{claim}[lemma]{Claim}
\newtheorem{proposition}[lemma]{Proposition}

\theoremstyle{definition}

\newtheorem{remark}[lemma]{Remark}

\begin{document}

\title{Metric structures in $L_1$:\\ Dimension, snowflakes, and average distortion}
\author{
    James R. Lee\thanks{Work partially supported by NSF grant CCR-0121555
    and an NSF Graduate Research Fellowship.  Part of this work was done while
    the author was an intern at Microsoft Research.}\\ U.C.
    Berkeley\\jrl@cs.berkeley.edu
    \and
    Manor Mendel\thanks{Work done while the author was a post-doc fellow at The
    Hebrew University, and supported in part by the Landau Center and by a grant
    from the Israeli Science Foundation (195/02).}\\University of
    Illinois\\mendelma@gmail.com
    \and
    Assaf Naor\\Microsoft Research\\anaor@microsoft.com
     }

\date{}

\maketitle

\begin{abstract}
We study the metric properties of finite subsets of $L_1$.  The
analysis of such metrics is central to a number of important
algorithmic problems involving the cut structure of weighted
graphs, including the Sparsest Cut Problem, one of the most
compelling open problems in the field of approximation algorithms.
Additionally, many open questions in geometric non-linear
functional analysis involve the properties of finite subsets of
$L_1$.
\smallskip

We present some new observations concerning the relation of $L_1$
to dimension, topology, and Euclidean distortion. We show that
every $n$-point subset of $L_1$ embeds into $L_2$ with average
distortion $O(\sqrt{\log n})$, yielding the first evidence that
the conjectured worst-case bound of $O(\sqrt{\log n})$ is valid.
We also address the issue of dimension reduction in $L_p$ for $p
\in (1,2)$.  We resolve a question left open in \cite{sahai} about
the impossibility of {\em linear} dimension reduction in the above
cases, and we show that the example of \cite{bc,ln} cannot be used
to prove a lower bound for the non-linear case. This is
accomplished by exhibiting constant-distortion embeddings of
snowflaked planar metrics into Euclidean space.
\end{abstract}

%\pagebreak

\section{Introduction}

This paper is devoted to the analysis of metric properties of
finite subsets of $L_1$. Such metrics occur in many important
algorithmic contexts, and their analysis is key to progress on
some fundamental problems.  For instance, an $O(\log
n)$-approximate max-flow/min-cut theorem proved elusive for many
years until, in \cite{LLR,AR}, it was shown to follow from a
theorem of Bourgain stating that every metric on $n$ points embeds
into $L_1$ with distortion $O(\log n)$.

 The importance of $L_1$
metrics has given rise to many problems and conjectures that have
attracted a lot of attention in recent years. to Four basic
problems of this type are as follows .
\begin{enumerate}[I{.}]
\item Is there an $L_1$ analog of the Johnson-Lindenstrauss
dimension reduction lemma~\cite{jl}? \item Are all $n$-point
subsets of $L_1$ $O\left(\sqrt{\log n}\right)$-embeddable into
Hilbert space? \item Are all squared-$\ell_2$ metrics
$O(1)$-embeddable into $L_1$? \item Are all planar graphs
$O(1)$-embeddable into $L_1$?
\end{enumerate}
(We recall that a squared-$\ell_2$ metric is a space $(X,d)$ for
which $(X, d^{1/2})$ embeds isometrically in a Hilbert space.)

Each of these questions has been asked many times before; we refer
to \cite{MatousekBook,HaifaOP,LinialICM,Indyk01}, in particular.
Despite an immense amount of interest and effort, the metric
properties of $L_1$ have proved quite elusive; hence the name
``The mysterious $L_1$'' appearing in a survey of Linial at the
ICM in 2002~\cite{LinialICM}.  In this paper, we attempt to offer
new insights into the above problems and touch on some
relationships between them. We refer the reader to the
book~\cite{MatousekBook} for an introductory account of the theory
of low distortion embeddings of metric spaces. In particular,
throughout this paper we shall use the standard terminology
appearing in~\cite{MatousekBook}.

\subsection{Results and techniques}

\noindent {\bf Euclidean distortion.}  Our first result addresses
problem (II) stated above.
We show that the answer to this question is positive on
average, in the following sense.

\begin{theorem}\label{thm:nL1} For every $f_1,\ldots,f_n\in
L_1$ there is a linear operator $T:L_1\to L_2$ such that
$$
\frac{\|T(f_i)-T(f_j)\|_2}{\|f_i-f_j\|_1} \geq \frac{1}{\sqrt{8 \log n}},
\qquad 1 \leq i < j \leq n, \textrm{ and } \\
$$
$$
\frac{1}{\binom{n}{2}} \sum_{1\le i<j\le
n}\left(\frac{\|T(f_i)-T(f_j)\|_2}{\|f_i-f_j\|_1}\right)^{1/2}\le
10.$$
\end{theorem}
In other words, for any $n$-point subset in $L_1$, there exists a
map into $L_2$ such that distances are contracted by at most
$O(\sqrt{\log n})$ and the average expansion is $O(1)$.
This yields the first positive evidence that
the conjectured worst-case bound of $O(\sqrt{\log n})$ holds.
%\begin{theorem}\label{thm:nL1}
%Let $X$ be an $n$-point subset of $L_1$. Then there exists a
%linear mapping $T:L_1\to L_2$ such that for every distinct $x,y\in
%X$, if we denote $D_{xy}=\frac{\|T(x)-T(y)\|_2}{\|x-y\|_1}$ then $
%D_{xy}\ge \Omega\left(\frac{1}{\sqrt{\log n}}\right)$ for every
%$x,y\in X$, and the average of the numbers
%$\left\{a_{xy}^{2/3}\right\}_{x,y\in X}$ is $O(1)$.
%\end{theorem}
We remark that a different notion of average embedding was
recently studied by Rabinovich~\cite{yuri}; there, one
tries to embed (planar) metrics into the line such that the {\em
average distance} does not change too much.

The exponent $1/2$ above has no significance, and we can actually
obtain the same result for any power $1-\e$, $\e>0$ (we refer to
Section~\ref{section:pisier} for details).
The proof of
Theorem~\ref{thm:nL1} follows from the following probabilistic
lemma, which is implicit in~\cite{mp}. We believe that this result
is of independent interest.
\begin{lemma} There exists a distribution over linear mappings $T:
L_1\to L_2$ such that for every $x\in L_1\setminus \{0\}$ the
random variable $\frac{\|T(x)\|_2}{\|x\|_1}$ has density
$\frac{e^{-1/(4x^2)}}{x^2\sqrt{\pi}}$.
\end{lemma}
In contrast to Theorem~\ref{thm:nL1}, we show that problem (II)
cannot be resolved positively using linear mappings. Specifically,
we show that there are arbitrarily large $n$-point subsets of
$L_1$ such that any linear embedding of them into $L_2$ incurs
distortion $\Omega(\sqrt{n})$. As a corollary we settle the
problem left open by Charikar and Sahai in~\cite{sahai}, whether
 dimension reduction with a {\em linear} map is possible in $L_p$, $p\notin
\{1,2\}$. The case $p=1$ was proved in~\cite{sahai} via linear
programming techniques, and it seems impossible to generalize
their method to arbitrary $L_p$. We show that there are
arbitrarily large $n$-point subsets  $X \subseteq L_p$ (namely,
the same point set used in~\cite{sahai} to handle the case $p=1$),
such that any linear embedding of $X$ into $\ell_p^d$ incurs
distortion $\Omega\left[(n/d)^{|1/p-1/2|}\right]$, thus dimension
reduction with a linear map is impossible in any $L_p$, $p\neq 2$.
Additionally, we show that there are arbitrarily large $n$-point
subsets  $X \subseteq L_1$ such any linear embedding of $X$ into
{\em any} $d$-dimensional normed space incurs distortion
$\Omega\left(\sqrt{n/d}\right)$. This generalizes the
Charikar-Sahai result to arbitrary low dimensional norms.
\bigskip

\noindent {\bf Dimension reduction.} In~\cite{bc}, and soon after
in~\cite{ln}, it was shown that if the Newman-Rabinovich diamond
graph on $n$ vertices $\alpha$-embeds into $\ell_1^d$ then $d\ge
n^{\Omega(1/\alpha^2)}$. The proof in~\cite{bc} is based on a
linear programming argument, while the proof in~\cite{ln} uses a
geometric argument which reduces the problem to bounding from
below the distortion required to embed the diamond graph in
$\ell_p$, $1<p<2$. These results settle the long standing open
problem of whether there is an $L_1$ analog of the
Johnson-Lindenstrauss dimension reduction lemma~\cite{jl}. (In
other words, they show that the answer to problem (I) above is
{\em No}.). In Section \ref{section:laakso}, we show that the
method of proof in~\cite{ln} can be used to provide an even more
striking counter example to this problem.

A metric space $X$ is called {\em doubling} with constant $C$ if
every ball in $X$ can be covered by $C$ balls of half the radius.
Doubling metrics with bounded doubling constants are widely viewed
as low dimensional (see \cite{GKL03,KLNN03} for some practical and
theoretical applications of this viewpoint).
%In fact, they
%have bounded Assouad dimension (see~\cite{Heinonen01} for the
%definition).
On the other hand, the doubling constant of the diamond
graphs is $\Omega(\sqrt{n})$ (where $n$ is the number of points).
Based on a
fractal construction due to Laakso~\cite{laakso} and the method
developed in~\cite{ln}, we prove the following theorem, which
shows a strong lower bound on the dimension required to represent
uniformly doubling subsets of $L_1$.

\begin{theorem}\label{thm:laakso} There are arbitrarily large
$n$-point subsets $X \subseteq L_1$ which are doubling with
constant $6$ but such that every $\alpha$-embedding of $X$ into
$\ell_1^d$ requires $d\ge n^{\Omega(1/\alpha^2)}$.
\end{theorem}

In~\cite{urs,GKL03} it was asked whether any subset of $\ell_2$
which is doubling well-embeds into $\ell_2^d$ (with bounds on the
distortion and the dimension that depend only on the doubling
constant).  In \cite{GKL03}, it was shown that a similar property
cannot hold for $\ell_1$.  Our lower bound exponentially
strengthens that result.

\bigskip

\noindent {\bf Planar metrics.} Our final result addresses
problems (III) and (IV).   Our motivation was an attempt to
generalize the argument in~\cite{ln} to prove that dimension
reduction is impossible in $L_p$ for any $1< p<2$. A natural
approach to this problem is to consider the point set used
in~\cite{bc,ln} (namely, a natural realization of the diamond
graph, $G$, in $L_1$) with the metric induced by the $L_p$ norm
instead of the $L_1$ norm. This is easily seen to amount to
proving lower bounds on the distortion required to embed the
metric space $(G,d_G^{1/p})$ in $\ell_p^h$. Unfortunately, this
approach cannot work since we show that, for any planar metric
$(X,d)$ and any $0<\e<1$, the metric space $(X,d^{1-\e})$ embeds
in Hilbert space with distortion $O\left(1/\sqrt{\e}\right)$, and
then using results of Johnson and Lindenstrauss \cite{jl}, and
Figiel, Lindenstrauss and Milman \cite{flm}, we conclude that this
metric can be $O(1/\sqrt{\e})$ embedded in $\ell_p^h$, where
$h=O(\log n)$. The proof of this interesting fact is a
straightforward application of Assouad's classical embedding
theorem~\cite{Assouad83} and Rao's embedding method~\cite{Rao99}.
The $O\left(1/\sqrt{\e}\right)$ upper bound is shown to be tight
for every value $0 < \varepsilon < 1$. We note that the case
$\e=1/2$ has been previously observed by A. Gupta in his
(unpublished) thesis~\cite{gupta-thesis}.

\remove{ \jnote{The following comment is only a little true.} It
follows that any planar metric embeds into squared-$\ell_2$ with
$O(1)$ distortion so that a positive solution to problem (III)
above implies a positive solution to problem (IV). }

\section{Average distortion Euclidean embedding of subsets of
$L_1$}\label{section:pisier}

The heart of our argument is the following lemma which is implicit
in~\cite{mp}, and which seems to be of independent interest.

\begin{lemma}\label{lem:representation} For every $0<p\le 2$ there is a probability space
$(\Omega,P)$ such that for every $\omega\in \Omega$ there is a
linear operator $T_\omega:L_p\to L_2$ such that for every $x\in
L_p\setminus \{0\}$ the random variable
$X=\frac{\|T_\omega(x)\|_2}{\|x\|_p}$ satisfies for every $a\in
\R$, $\E e^{-a X^2}=e^{-a^{p/2}}$. In particular, for $p=1$ the
density of $X$ is $\frac{e^{-1/(4x^2)}}{x^2\sqrt{\pi}}$.
\end{lemma}

\begin{proof} Consider the following three
sequences of random variables, $\{Y_j\}_{j\ge 1}$,
$\{\theta_j\}_{j\ge 1}$, $\{g_j\}_{j\ge 1}$, such that each
variable is independent of the others. For each $j\ge 1$, $Y_j$ is
uniformly distributed on $[0,1]$, $g_j$ is a standard Gaussian and
$\theta_j$ is an exponential random variable, i.e. for $\lambda\ge
0$, $P(\theta_j>\lambda)=e^{-\lambda}$. Set
$\Gamma_j=\theta_1+\cdots+\theta_j$. By Proposition 1.5. in
\cite{mp}, there is a constant $C=C(p)$ such that if we define for
$f\in L_p$
$$
V(f)=C\sum_{j\ge 1} \frac{g_j}{\Gamma_j^{1/p}}f(Y_j),
$$
then $\E e^{iV(f)}=e^{-\|f\|_p^p}$.

Assume that the random variables $\{Y_j\}_{j\ge 1}$ and
$\{\Gamma_j\}_{j\ge 1}$ are defined on a probability space
$(\Omega,P)$ and that $\{g_j\}_{j\ge 1}$ are defined on a
probability space $(\Omega',P')$, in which case we use the
notation $V(f)=V(f;\omega,\omega')$. Define for $\omega\in \Omega$
a linear operator $T_\omega: L_p\to L_2(\Omega',P')$ by
$T_\omega(f)= V(f;\omega,\cdot)$. Since for every fixed $\omega\in
\Omega$ the random variable $V(f;\omega,\cdot)$ is Gaussian with
variance $\|T_\omega(f)\|_2^2$, for every $a\in \R$, $\E_{P'}
e^{iaV(s;\omega,\cdot)}= e^{-a^2\|T_\omega(f)\|_2^2}$. Taking
expectation with respect to $P$ we find that, $\E_P
e^{-a^2\|T_\omega(f)\|_2^2}=e^{-a^p\|f\|_p^p}$. This implies the
required identity. The explicit distribution in the case $p=1$
follows from the fact that the inverse Laplace transform of
$x\mapsto e^{-\sqrt{x}}$ is $y\mapsto
\frac{e^{-1/(4y)}}{2\sqrt{\pi y^3}}$ (see for
example~\cite{widder,durrett}).
\end{proof}

\begin{comment}

We now pass to case case $p=1$, although all our argument carry
over to the case of general $1<p<2$

Using the notation of Lemma~\ref{lem:representation}, a standard
application of Markov's inequality shows that for every $\e>0$,
$P\{X\le \e\}\le e^{-c/\e^\alpha}$ where $c$ depends only on $p$
and $\alpha=\frac{2p}{2-p}$. It follows that for every
$f_1,\ldots, f_n\in L_p$, with probability greater than $1/2$, for
every $1\le i<j\le n$, $\|f_i-f_j\|_p\le O\left[(\log
n)^{1/p-1/2}\right]\cdot \|T_\omega(f_i)-T_\omega(f_j)\|_2$. On
the other hand, since for $1<p<2$ $\E X<\infty$ and for every
$1\le i<j\le n$, $\E\|T_\omega(f_i)-T_\omega(f_j)\|_2=(\E
X)\|f_i-f_j\|_p$, a simple averaging argument proves the estimate
$e_n(f_i,f_j)=C(p)(\log n)^{1/p-1/2}$ by reducing it to the
Hilbertian case (i.e. Kirszbraun's theorem). We refer to Theorem
2.12 in~\cite{mp} for the details. We record below another
application of Lemma~\ref{lem:representation}. It has long been
conjectured that any $n$ point subset of $L_1$ embeds in Hilbert
space with distortion $O\left(\sqrt{\log n}\right)$. The following
corollary shows that this is indeed the case when we only require
a bound on the average distortion:
\end{comment}

\begin{proof}[Proof of Theorem~\ref{thm:nL1}] Using the notation of
lemma~\ref{lem:representation} (in the case $p=1$) we find that
for every $a>0$, $\E e^{-a X^2}=e^{-\sqrt{a}}$. Hence, for every
$a,\e>0$ and every $1<i<j\le n$,
\begin{eqnarray*}
P\left(\frac{\|T_\omega(f_i)-T_\omega(f_j)\|_2}{\|f_i-f_j\|_1}\le
\e\right)=P\left(e^{-a X^2}\ge e^{-a\e^2}\right)\le
e^{a\e^2-\sqrt{a}}.
\end{eqnarray*}
Choosing $a=\frac{1}{4\e^4}$ the above upper bound becomes
$e^{-1/(4\e^2)}$. Consider the set
$$
A=\bigcap_{1\le i<j\le n}
\left\{\frac{\|T_\omega(f_i)-T_\omega(f_j)\|_2}{\|f_i-f_j\|_1}\ge
\frac{1}{\sqrt{8\log n}}\right\}\subseteq\Omega.
$$
By the union bound, $P(A)>\frac12$, so that
\begin{multline*}
\frac{1}{P(A)}\E\left[\frac{1}{\binom{n}{2}} \sum_{1\le i<j\le
n}\left(\frac{\|T_\omega(f_i)-T_\omega(f_j)\|_2}{\|f_i-f_j\|_1}\right)^{1/2}\right]
\le 2\E X^{1/2}=\frac{2}{\sqrt{\pi}}\int_0^\infty x^{1/2}\cdot
\frac{e^{-1/(4x^2)}}{x^2}dx<10.
\end{multline*}
It follows that there exists $\omega\in A$ for which the operator
$T=T_\omega$ has the desired properties.
\end{proof}

\begin{remark} There is nothing special about the choice of
the power $1/2$ in Theorem~\ref{thm:nL1}. When $p=1$, $\E
X=\infty$ but $\E X^{1-\e}<\infty$ for every $0<\e<1$, so we may
write the above average with the power $1-\e$ replacing the
exponent $1/2$. Obvious generalizations of Theorem~\ref{thm:nL1}
hold true for every $1<p<2$, in which case the average distortion
is of order $C(p)(\log n)^{1/p-1/2}$ (and the power can be taken
to be $1$).
\end{remark}

\section{The impossibility of dimension reduction with a {\em linear} map in
$L_p$, $p\neq 2$}

The above method cannot yield a $O\left(\sqrt{\log n}\right)$
bound on the Euclidean distortion of $n$-point subsets of $L_1$.
In fact, there are arbitrarily large $n$-point subsets of $L_1$ on
which any {\em linear} embedding into $L_2$ incurs distortion at
least $\sqrt{\frac{n-1}{2}}$. This follows from the following
simple lemma:

\begin{lemma}\label{lem:nolin} For every $1\le p\le\infty$ there
are arbitrarily large $n$-point subsets of $L_p$ on which any
linear embedding into $L_2$ incurs distortion at least
$\left(\frac{n-1}{2}\right)^{\left|1/p-1/2\right|}$.
\end{lemma}

\begin{proof} Let $w_1,\ldots,w_{2^k}$ be
the rows of the $2^k\times 2^k$ Walsh matrix (i.e. the simplest
Hadamard matrix). Write $w_i=\sum_{j=1}^{2^k} w_{ij} e_j$ where
$e_1,\ldots, e_{2^k}$ are the standard unit vectors in $\R^{2^k}$.
Consider the set
$A=\{0\}\cup\{w_i\}_{i=1}^{2^k}\cup\{e_i\}_{i=1}^{2^k}\subset
\ell_p$. Let $T:\ell_p\to L_2$ be any linear operator which is non
contracting and $L$-Lipschitz on $A$. Assume first of all that
$1\le p<2$. Then:
\begin{multline*}
2^{k(1+2/p)}= \sum_{i=1}^{2^k}\| w_i\|_p^2 \le \sum_{i=1}^{2^k}\|T
w_i\|_2^2= \sum_{i=1}^{2^k}\left\|\sum_{j=1}^{2^k}
w_{ij}T(e_j)\right\|_2^2\\
= \sum_{i=1}^{2^k}\sum_{j=1}^{2^k}\left\langle
w_i,w_j\right\rangle \left\langle
T(e_i),T(e_j)\right\rangle=2^k\sum_{j=1}^{2^k} \|T(e_j)\|_2^2\le
4^k\cdot L^2,
\end{multline*}
which implies that $L\ge
2^{k(1/p-1/2)}=\left(\frac{|A|-1}{2}\right)^{1/p-1/2}$. When $p>2$
apply the same reasoning, with the inequalities reversed.
\end{proof}

We remark that the above point set was also used by Charikar and
Sahai~\cite{sahai} to give a lower bound on dimension reduction
with a linear map in $L_1$. Their proof used a linear programming
argument, which doesn't seem to be generalizable to the the case
of $L_p$, $p>1$. Lemma~\ref{lem:nolin} formally implies their
result (with a significantly simpler proof), and in fact proves
the impossibility of dimension reduction with a linear map in any
$L_p$, $p\neq 2$. Indeed, if there were a linear operator which
embeds $A$ into $\ell_p^d$ with distortion $D$ then it would also
be a $D\cdot d^{|1/p-1/2|}$ embedding into $\ell_2^d$. It follows
that $D\ge \left(\frac{|A|-1}{2d}\right)^{|1/p-1/2|}$. Similarly,
since by John's theorem (see e.g.~\cite{mil-schecht}) any
$d$-dimensional normed space is $\sqrt{d}$ equivalent to Hilbert
space, we deduce that there are arbitrarily large $n$-point
subsets of $L_1$, any linear embedding of which into any
$d$-dimensional normed space incurs distortion at least
$\sqrt{\frac{n-1}{2d}}$.

\section{An inherently high-dimensional doubling metric in
$L_1$}\label{section:laakso}

This section is devoted to the proof of Theorem~\ref{thm:laakso}.

\begin{proof}[Proof of Theorem~\ref{thm:laakso}]
Consider the Laakso graphs, $\{G_i\}_{i=0}^\infty$, which are
defined as follows. $G_0$ is the graph on two vertices with one
edge. To construct $G_i$, take six copies of $G_{i-1}$ and scale
their metric by a factor of $\frac14$. We glue four of them
cyclicly by identifying pairs of endpoints, and attach at two
opposite gluing points the remaining two copies. See Figure 1
below.

\begin{figure}[ht]
\bigskip
\ \centering
\input{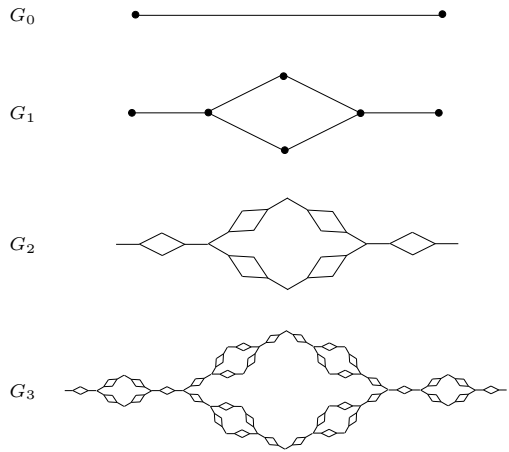}
\caption{The Laakso graphs.} \label{fig:lang}
\end{figure}

As shown in~\cite{laakso}, the graphs $\{G_i\}_{i=0}^\infty$ are
uniformly doubling (see also~\cite{urs}, for a simple argument
showing they are doubling with constant $6$). Moreover, since the
$G_i$'s are series parallel graphs, they embed uniformly in $L_1$
(see~\cite{sinclair}).

We will show below that any embedding of $G_i$ in $L_p$, $1<p\le
2$ incurs distortion at least $\sqrt{1+\frac{p-1}{4}i}$. We then
conclude as in~\cite{ln} by observing that $\ell_1^d$ is
$3$-isomorphic to $\ell_p^d$ when $p=1+\frac{1}{\log d}$, so that
if $G_i$ embeds with distortion $\alpha$ in $\ell_1^d$ then
$\alpha\ge \sqrt{\frac{i}{40\log d}}$. This implies the required
result since $i\approx \log |G_i|$.

The proof of the lower bound for the distortion required to embed
$G_i$ into $L_p$ is by induction on $i$. We shall prove by
induction that whenever $f:G_i\to L_p$ is non-contracting then
there exist two adjacent vertices $u,v\in G_i$ such that
$\|f(u)-f(v)\|_p\ge d_{G_i}(u,v)\sqrt{1+\frac{p-1}{4}i}$ (observe
that for $u,v\in G_{i-1}$, $d_{G_{i-1}}(u,v)=d_{G_i}(u,v)$). For
$i=0$ there is nothing to prove. For $i\ge 1$, since $G_i$
contains an isometric copy of $G_{i-1}$, there are $u,v\in G_i$
corresponding to two adjacent vertices in $G_{i-1}$ such that
$\|f(u)-f(v)\|_p\ge d_{G_i}(u,v)\sqrt{1+\frac{p-1}{4}(i-1)}$. Let
$a,b$ be the two midpoints between $u$ and $v$ in $G_i$. By Lemma
2.1 in~\cite{ln},
\begin{multline*}
\|f(u)-f(v)\|_p^2 + (p-1)\|f(a)-f(b)\|_p^2\\ \le
\|f(u)-f(a)\|_p^2+\|f(a)-f(v)\|_p^2+\|f(v)-f(b)\|_p^2+\|f(b)-f(u)\|_p^2.
\end{multline*}
Hence:
\begin{eqnarray*}
&&\!\!\!\!\!\!\!\!\!\!\!\!\!\max\{\|f(u)-f(a)\|_p^2,\|f(a)-f(v)\|_p^2,\|f(v)-f(b)\|_p^2,\|f(b)-f(u)\|_p^2\}\\&\ge&
\frac14\|f(u)-f(v)\|_p^2+\frac14(p-1)\|f(a)-f(b)\|_p^2\\&\ge&\frac14\left(1+\frac{p-1}{4}(i-1)\right)d_{G_i}(u,v)^2+
\frac{p-1}{4}d_{G_i}(a,b)^2\\
&=&\frac14\left(1+\frac{p-1}{4}i\right)d_{G_i}(u,v)^2\\&=&
 \left(1+\frac{p-1}{4}i\right)\max\{d_{G_i}(u,a)^2,d_{G_i}(a,v)^2,d_{G_i}(v,b)^2,d_{G_i}(b,u)^2\}.
\end{eqnarray*}

\end{proof}

We end this section by observing that the above approach also
gives a lower bound on the dimension required to embed expanders
in $\ell_\infty$.

\begin{proposition}\label{prop:dimexp} Let $G$ be an $n$-point constant degree expander
which embeds in
$\ell_\infty^d$ with distortion at most $\alpha$. Then $d\ge
n^{\Omega(1/\alpha)}$.
\end{proposition}

\begin{proof} By Matou\v{s}ek's lower bound for the distortion required to embed expanders in $\ell_p$~\cite{matexpander},
any embedding of $G$ into $\ell_p$ incurs distortion
$\Omega\left(\frac{\log n}{p}\right)$. Since $\ell_\infty^d$ is
$O(1)$-equivalent to $\ell_{\log d}^d$, we deduce that $\alpha\ge
\Omega\left(\frac{\log n}{\log d}\right)$.
\end{proof}

We can also obtain a lower bound on the dimension required to
embed the Hamming cube $\{0,1\}^k$ into $\ell_\infty$. Our proof
uses a simple concentration argument. An analogous concentration
argument yields an alternative proof of
Proposition~\ref{prop:dimexp}.

\begin{proposition}\label{prop:dimcube} Assume that $\{0,1\}^k$
embeds into $\ell_\infty^d$ with distortion $\alpha$. Then $d\ge
2^{k\Omega(1/\alpha^2)}$.
\end{proposition}

\begin{proof} Let $f=(f_1,\ldots,f_d):\{0,1\}^k\to \ell_\infty^d$ be a
contraction such that for every $u,v\in \{0,1\}^d$,
$\|f(u)-f(v)\|_\infty\ge \frac{1}{\alpha}d(u,v)$ (where
$d(\cdot,\cdot)$ denotes the Hamming metric). Denote by $P$ the
uniform probability measure on $\{0,1\}^k$. Since for every $1\le
i\le k$, $f_i$ is $1$-Lipschitz, the standard concentration
inequality on the hypercube (see~\cite{MatousekBook}) implies that
$P\left(|f_i(u)-\E f_i|\ge k/(4\alpha)\right)\le
2e^{-k/(32\alpha^2)}$. On the other hand, if $u,v\in \{0,1\}^k$
are such that $d(u,v)=k$ then there exist $1\le i\le d$ for which
$|f_i(u)-f_i(v)|\ge k/\alpha$, implying that $\max\{|f_i(u)-\E
f_i|,|f_i(v)-\E f_i|\}>k/(4\alpha)$. By the union bound it follows
that $de^{-\Omega(k/\alpha^2)}\ge 1$, as required.
\end{proof}

\section{Snowflake versions of planar metrics}

The problem of whether there is an analog of the
Johnson-Lindenstrauss dimension reduction lemma in $L_p$, $1<p<2$,
is an interesting one which remains open. In view of the above
proof and the proof in~\cite{ln}, a natural point set which is a
candidate to demonstrate the impossibility of dimension reduction
in $L_p$ is the realization of the diamond graph in $\ell_1$ which
appears in~\cite{bc}, equipped with the $\ell_p$ metric. Since
this point set consists of vectors whose coordinates are either 0 or 1 (i.e.
subsets of the cube), this amounts to
considering the diamond graph with its metric raised to the power
$\frac1p$. Unfortunately, this approach cannot work; we show below
that any planar graph whose metric is raised to the power $1-\e$
has Euclidean distortion $O\left(1/\sqrt{\e}\right)$.

Given a metric space $(X,d)$ and $\e>0$, the metric space
$(X,d^{1-\e})$ is known in geometric analysis (see
e.g.~\cite{Heinonen01}) as the $1-\e$ snowflake version of
$(X,d)$. Assouad's classical theorem~\cite{Assouad83} states that
any snowflake version of a doubling metric space is bi-Lipschitz
equivalent to a subset of some finite dimensional Euclidean space.
A quantitative version of this result (with bounds on the
distortion and the dimension) was obtained in~\cite{GKL03}. The
following theorem is proved by combining embedding techniques of
Rao~\cite{Rao99} and Assouad~\cite{Assouad83}.  A similar analysis
is also used in~\cite{GKL03}. In what follows we call a metric
$K_r$-excluded if it is the metric on a subset of a weighted graph
which does not admit a $K_r$ minor. In particular, planar metrics
are all $K_5$-excluded.

\begin{theorem}\label{thm:snowflake} For any
$r\in\mathbb{N}$ there exists a constant $C(r)$ such that for
every $0 < \epsilon < 1$, a $1-\e$ snowflake version of a
$K_r$-excluded metric embeds into $\ell_2$ with distortion at most
$C(r)/\sqrt{\e}$.
\end{theorem}

Our argument is based on the following lemma, the proof of which
is contained in~\cite{Rao99}.

\begin{lemma}\label{lem:rao} For every $r\in \mathbb{N}$ there is
a constant $\delta=\delta(r)$ such that for every $\rho>0$ and
every $K_r$-excluded metric $(X,d)$ there exists a finitely
supported probability distribution $\mu$ on partitions of $X$ with
the following properties:
\begin{enumerate}
\item For every $P \in \supp(\mu)$, and for every $C \in P$,
$\diam(C) \leq \rho$. \item For every $x\in X$, $\E_{\mu}
\sum_{C\in P} d(x,X\setminus C) \geq \delta \rho$.
\end{enumerate}
\end{lemma}

Observe that the sum under the expectation in (2) above actually
consists of only one summand.

\begin{proof}[Proof of Theorem~\ref{thm:snowflake}] Let $X$ be a
$K_r$-excluded metric. For each $n \in \mathbb Z$, we define a map
$\phi_n$ as follows. Let $\mu_n$ be the probability distribution
on partitions of $X$ from Lemma~\ref{lem:rao} with
$\rho=2^{n/(1-\e)}$. Fix a partition $P \in \supp(\mu_n)$. For any
$\sigma \in \{-1,+1\}^{|P|}$, consider $\sigma$ to be indexed by
$C \in P$ so that $\sigma_C$ denotes the value of $\sigma$ at $C$.
Following Rao~\cite{Rao99}, define
$$
\phi_P(x) = \bigoplus_{\sigma \in \{-1,+1\}^{|P|}}
\sqrt{\frac{1}{2^{|P|}}} \sum_{C \in P} \sigma_C \cdot
d(x,X\setminus C),
$$
and write $\phi_n = \bigoplus_{P \in \supp(\mu_n)}
\sqrt{\mu_n(P)}\, \phi_P$ (here the symbol $\oplus$ refers to the
concatenation operator).

Now, following Assouad~\cite{Assouad83}, let $\{e_i\}_{i\in
\mathbb{Z}}$ be an orthonormal basis of $\ell_2$, and set
$$
\Phi(x) = \sum_{n \in \mathbb Z} 2^{-n\e/(1-\e)} \phi_n(x) \otimes
e_n$$

\begin{claim}
For every $n \in \mathbb Z$, and $x,y \in X$, we have $||\phi_n(x)
- \phi_n(y)||_2 \leq 2\cdot\min\left\{d(x,y),
2^{n/(1-\e)}\right\}$. Additionally, if $d(x,y)
> 2^{n/(1-\e)}$, then $||\phi_n(x) - \phi_n(y)||_2 \geq \delta\, 2^{n/(1-\e)}$.
\end{claim}

\begin{proof}
For any partition $P \in \supp(\mu_n)$, let $C_x, C_y$ be the
clusters of $P$ containing $x$ and $y$, respectively.  Note that
since for every $C \in P$, $\diam(C) \leq 2^{n/(1-\e)}$, when
$d(x,y)
> 2^{n/(1-\e)}$, we have $C_x \neq C_y$.  In this case,
\begin{eqnarray*}||\phi_P(x) - \phi_P(y)||_2^2 &=& \E_{\sigma \in
\{-1,+1\}^{|P|}} |\sigma_{C_x} d(x,X\setminus C_x) - \sigma_{C_y}
d(y,X\setminus C_y)|^2 \\&\geq& \frac{d(x,X\setminus
C_x)^2+d(y,X\setminus C_y)^2}{2}.
\end{eqnarray*} It follows that
\begin{eqnarray*}
||\phi_n(x) - \phi_n(y)||_2^2 &=& \E_{\mu_n} ||\phi_P(x) -
\phi_P(y)||_2^2\\& \geq& \frac{\E_{\mu_n} d(x,X\setminus
C_x)^2+\E_{\mu_n}d(y,X\setminus C_y)^2}{2} \geq \left(\delta
\,2^{n/(1-\e)}\right)^2.
\end{eqnarray*}

 On the other hand, for
every $x,y \in X$, since $d(x,X\setminus C_x), d(y,X\setminus C_y)
\leq 2^{n/(1-\e)}$, we have that
 $||\phi_P(x) - \phi_P(y)||_2 \leq 2\cdot
\min\left\{d(x,y), 2^{n/(1-\e)}\right\}$, hence $||\phi_n(x) -
\phi_n(y)||_2 \leq 2 \cdot \min\left\{d(x,y),
2^{n/(1-\e)}\right\}$.
\end{proof}

To finish the analysis, let us fix $x,y \in X$ and let $m$ be such
that $d(x,y)^{1-\e} \in \left(2^{m}, 2^{m+1}\right]$. In this
case,

\begin{comment}
\begin{eqnarray*}
||\Phi(x) - \Phi(y)||_2^2& = & \sum_{i=0}^{k-1}\left\|\sum_{n \in
\mathbb Z} 2^{-\epsilon(nk+i)} \left(\phi_{nk+i}(x) -
\phi_{nk+i}(y)\right)\right\|_2^2\\
&\le& \sum_{i=0}^{k-1}\left(\sum_{n\in \mathbb{Z}}
2^{-\epsilon(nk+i)} \|\phi_{nk+i}(x) -
\phi_{nk+i}(y)\|_2\right)^2\\&\leq & \sum_{i=0}^{k-1} 4\left(
\sum_{n:\ nk + i < m} 2^{(1-\epsilon) (nk+i)} + d(x,y) \sum_{n :\
nk + i \geq m}
2^{-\epsilon(nk+i)}\right)^2 \\
%4^{(\epsilon-1)n} ||\phi_n(x) - \phi_n(y)||_2 \\
&\leq & 4\sum_{i=0}^k\left(\frac{2^{(1-\e)m}}{1-2^{-(1-\e)k}}+ d(x,y)\frac{2^{-\e m}}{1-2^{-\e k}} \right)^2 \\
&= &  O\left[\frac{\log(1/\delta)}{\e(1-\e)}\right]\cdot
d(x,y)^{2(1-\e)}
% 4k\left(\frac{1}{1-4^{-\epsilon k}}
%4^{\epsilon m} + \frac{1}{1-4^{-\epsilon k}}\, d(x,y)^2
%4^{(\epsilon-1)m}\right)
\end{eqnarray*}
\end{comment}

\begin{eqnarray*}
||\Phi(x) - \Phi(y)||_2^2& = & \sum_{n \in \mathbb
Z}2^{-2n\e/(1-\e)} \left \| \phi_{n}(x) -
\phi_{n}(y)\right\|_2^2\\
&\leq &  4\sum_{ n < m} 2^{2 n} + 4d(x,y)^2 \sum_{ n \geq m}
2^{-2n\epsilon/(1-\e)}  \\
&= & 2^{2m+1}+ 4d(x,y)^2\frac{2^{-2m\e/(1-\e)}}{1-2^{-2\e/(1-\e) }}  \\
&= &  O\left(1/\e\right)\cdot d(x,y)^{2(1-\e)}
% 4k\left(\frac{1}{1-4^{-\epsilon k}}
%4^{\epsilon m} + \frac{1}{1-4^{-\epsilon k}}\, d(x,y)^2
%4^{(\epsilon-1)m}\right)
\end{eqnarray*}

On the other hand,
\begin{eqnarray*}
\|\Phi(x) - \Phi(y)||_2 \geq   2^{-m\epsilon/(1-\e)} \|\phi_{m}(x)
-\phi_{m}(y)\|_2  \ge \delta 2^m\ge \frac{\delta}{2}d(x,y)^{1-\e}.
\end{eqnarray*}
The proof is complete.
\end{proof}

\begin{remark} The $O\left(1/\sqrt{\e}\right)$ upper bound in
Theorem~\ref{thm:snowflake} is tight. In fact, for $i\approx
1/\e$, the $1-\e$ snowflake version of the Laakso graph $G_i$
(presented in Section~\ref{section:laakso}) has Euclidean
distortion $\Omega\left(1/\sqrt{\e}\right)$. To see this, let
$f:G_i\to \ell_2$ be any non-contracting embedding of
$(G_i,d_{G_i}^{1-\e})$ into $\ell_2$. For $j\le i$ denote by $K_j$
the Lipschitz constant of the restriction of $f$ to $(G_j,
d_{G_i}^{1-\e})$ (as before, we think of $G_j$ as a subset of
$G_i$). Clearly $K_0=1$, and the same reasoning as in the proof of
Theorem~\ref{thm:laakso} shows that for $j\ge 1$, $K_j^2\ge
\frac{K_{j-1}^2}{4^\e}+\frac14$. This implies that $K_i^2\ge
\frac14+\frac{1}{4^\e}+\ldots+\frac{1}{4^{i\e}}=\Omega(1/\e)$, as
required.

\end{remark}

\bibliographystyle{abbrv}

\begin{thebibliography}{10}

\bibitem{Assouad83}
P.~Assouad.
\newblock Plongements lipschitziens dans {${\bf R}\sp{n}$}.
\newblock {\em Bull. Soc. Math. France}, 111(4):429--448, 1983.

\bibitem{AR}
Y.~Aumann and Y.~Rabani.
\newblock An ${O}(\log k)$ approximate min-cut max-flow theorem and
  approximation algorithm.
\newblock {\em SIAM J. Comput.}, 27(1):291--301, 1998.

\bibitem{bc}
M.~Charikar and B.~Brinkman.
\newblock On the impossibility of dimension reduction in $\ell_1$.
\newblock In {\em to appear in Proceedings of the 44th Annual IEEE Conference
  on Foundations of Computer Science}, pages 514--523. ACM, 2003.

\bibitem{sahai}
M.~Charikar and A.~Sahai.
\newblock Dimension reduction in the $\ell_1$ norm.
\newblock In {\em Proceedings of the 43rd Annual IEEE Conference on Foundations
  of Computer Science}, pages 251--260. ACM, 2002.

\bibitem{durrett}R. Durrett.
\newblock Probability: theory and examples. Second edition.
\newblock {\em Duxbury Press, Belmont, CA}, 1996

\bibitem{flm}
T.~Figiel, J.~Lindenstrauss, and V.~D. Milman.
\newblock The dimension of almost spherical sections of convex bodies.
\newblock {\em Acta Math.}, 139(1-2):53--94, 1977.

\bibitem{gupta-thesis} A.~Gupta.
\newblock Embeddings of Finite Metrics.
\newblock {\em Ph.D. thesis, University of California, Berkeley},
2000.


\bibitem{sinclair}
A.~Gupta, I.~Newman, Y.~Rabinovich, and A.~Sinclair.
\newblock Cuts, trees and $\ell_1$ embeddings.
\newblock In {\em Proceedings of the 40th Annual Symposium on Foundations of
  Computer Science}, pages 399--409, 1999.

\bibitem{GKL03}
A.~Gutpa, R.~Krauthgamer, and J.~R. Lee.
\newblock Bounded geometries, fractals, and low-distortion embeddings.
\newblock In {\em Proceedings of the 44th Annual Symposium on Foundations of
  Computer Science}, pages 534--543, 2003.

\bibitem{Heinonen01}
J.~Heinonen.
\newblock {\em Lectures on analysis on metric spaces}.
\newblock Universitext. Springer-Verlag, New York, 2001.

\bibitem{Indyk01}
P.~Indyk.
\newblock Algorithmic applications of low-distortion geometric embeddings.
\newblock In {\em Proceedings of the 42nd Annual IEEE Symposium on Foundations
  of Computer Science}, pages 10--33. October 2001.

\bibitem{jl}
W.~B. Johnson and J.~Lindenstrauss.
\newblock Extensions of {L}ipschitz mappings into a {H}ilbert space.
\newblock In {\em Conference in modern analysis and probability (New Haven,
  Conn., 1982)}, volume~26 of {\em Contemp. Math.}, pages 189--206. Amer. Math.
  Soc., Providence, RI, 1984.

\bibitem{KLNN03}
R.~Krauthgamer and J.~R. Lee.
\newblock Navigating nets: Simple algorithms for proximity search.
\newblock Submitted, 2003.

\bibitem{laakso}
T.~J. Laakso.
\newblock Ahlfors {$Q$}-regular spaces with arbitrary {$Q>1$} admitting weak
  {P}oincar\'e inequality.
\newblock {\em Geom. Funct. Anal.}, 10(1):111--123, 2000.

\bibitem{urs}
U.~Lang and C.~Plaut.
\newblock Bilipschitz embeddings of metric spaces into space forms.
\newblock {\em Geom. Dedicata}, 87(1-3):285--307, 2001.

\bibitem{ln}
J.~R. Lee and A.~Naor.
\newblock Embedding the diamond graph in ${L}_p$ and dimension reduction in
  ${L}_1$.
\newblock {\em To appear in Geometric and Functional Analysis}, 2003.

\bibitem{LinialICM}
N.~Linial.
\newblock Finite metric spaces - combinatorics, geometry and algorithms.
\newblock In {\em Proceedings of the International Congress of Mathematicians
  III}, pages 573--586, 2002.

\bibitem{LLR}
N.~Linial, E.~London, and Y.~Rabinovich.
\newblock The geometry of graphs and some of its algorithmic applications.
\newblock {\em Combinatorica}, 15(2):215--245, 1995.

\bibitem{mp}
M.~B. Marcus and G.~Pisier.
\newblock Characterizations of almost surely continuous {$p$}-stable random
  {F}ourier series and strongly stationary processes.
\newblock {\em Acta Math.}, 152(3-4):245--301, 1984.

\bibitem{matexpander}
J.~Matou{\v{s}}ek.
\newblock On embedding expanders into {$l\sb p$} spaces.
\newblock {\em Israel J. Math.}, 102:189--197, 1997.

\bibitem{MatousekBook}
J.~Matou{\v{s}}ek.
\newblock {\em Lectures on discrete geometry}, volume 212 of {\em Graduate
  Texts in Mathematics}.
\newblock Springer-Verlag, New York, 2002.

\bibitem{HaifaOP}
J.~Matou{\v{s}}ek.
\newblock Open problems, workshop on discrete metric spaces and their
  algorithmic appl ications.
\newblock Haifa, March 2002.

\bibitem{mil-schecht}
V.~D. Milman and G.~Schechtman.
\newblock {\em Asymptotic theory of finite-dimensional normed spaces}.
\newblock Springer-Verlag, Berlin, 1986.
\newblock With an appendix by M. Gromov.

\bibitem{yuri}
Y.~Rabinovich.
\newblock On average distorsion of embedding metrics into the line and into
  $l_1$.
\newblock In {\em Proceedings of the 35th Annual ACM Symposium on Theory of
  Computing}, pages 465--462. ACM, 2003.

\bibitem{Rao99}
S.~Rao.
\newblock Small distortion and volume preserving embeddings for planar and
  {E}uclidean metrics.
\newblock In {\em Proceedings of the 15th Annual Symposium on Computational
  Geometry}, pages 300--306. ACM, 1999.

\bibitem{widder}
D.~V. Widder.
\newblock {\em The {L}aplace {T}ransform}.
\newblock Princeton Mathematical Series, v. 6. Princeton University Press,
  Princeton, N. J., 1941.

\end{thebibliography}

\def\cprime{$'$} \def\cprime{$'$} \def\cprime{$'$}

\end{document}